# SYMMETRY OF THE TRIPLE OCTONIONIC PRODUCT

**© 2018 г. M. V. Kharinov**
*(199178, Russia, St. Petersburg,*
*14 line of V.I., 39, SPIIRAS)*
*e-mail: khar@iias.spb.su*

*The Hermitian decomposition of a linear operator is generalized to the case of two or more operations. An additive expansion of the product of three octonions into three parts is constructed, wherein each part either preserve or change the sign under the action of the Hermitian conjugation and operation of inversion of the multiplicative order of three hypercomplex numbers, as well as under the composition of specified operations. The product of three octonions, in particular quaternions, with conjugate central factor is presented as the sum of mutually orthogonal triple anticommutator, triple commutator and associator that vanishes in the case of associative quaternions. The triple commutator is treated as a generalization of the cross product to the case of three arguments both for quaternions and octonions. A generalized cross product is introduced as an antisymmetric component of the triple octonionic product that changes sign both for inversion of the multiplicative order of three arguments, and for the Hermitian conjugation of the product considered respectively to the central argument. The definition of the cross product of three hypercomplex numbers deduced from symmetry considerations is compared with the solution of S. Okubo (1993) and the modern solution of T. Dray and M.A. Corinne (2015). It is shown that the derived definition is equivalent to the first solution presented by S.Okubo in an insufficiently perfect form.*

## 1. INTRODUCTION

In this paper the generalization of a cross product to the case of three arguments for a four-dimensional and eight-dimensional vector spaces is considered. At present, one of the obstacles to solving this kind of problem is the commonly used definition of a vector product by means of an intuitively percepted "right-hand rule". In the toolkit of *hypercomplex numbers* (non-commutative quaternions and non-associative octonions [1, 2]), the notion of a cross product of a pair of vectors is introduced without referring to intuition, which simplifies the solving of the problem. However, unlike the cross product itself, acquaintance of the reader with hypercomplex numbers is usually limited, since the hypercomplex numbers are mentioned only in fewer courses on abstract algebra [3]. To quickly master the hypercomplex numbers it is enough to use the popular book [2]. Many useful rules for working with hypercomplex numbers one may find in [4]. In the presented paper only formulas that most important for a context are listed. Unnumbered formulas are mainly given either for memorization, or for explaining the meaning of the notation.

In the paper, we adhere to the notation of [2]. The aim of the paper is to generalize the additive Hermitian decomposition, to substantiate the generalization of the cross product, and to simplify the toolbox of hypercomplex numbers, which can prove to be especially useful when working with non-associative octonions.

## 2. ELEMENTARY INFORMATION

By $i_0$ we denote the multiplicative unit that commutes with any hypercomplex number $u$ and, when multiplied, leaves it unchanged:

$$i_0 u = u i_0 = u .$$

Note that $i_0$ is a simple renaming of 1 and, say $5i_0 \equiv 5\cdot 1 \equiv 5$.

Let $\bar{u}$ denotes the hypercomplex number which is conjugate to the number $u$, and is related to $u$ by the formula [4]:

$$\bar{u} = 2(u, i_0)i_0 - u ,$$

where $(u, i_0)$ is the inner product of the vector $u$ and vector $i_0$.

In four algebras of *hypercomplex numbers* (real numbers, complex numbers, quaternions and octonions), the square of length $(u,u)$ of the vector $u$ is introduced as the product of the vector $u$ and the conjugate vector $\bar{u}$, and the inner product of the vectors $u_1$ and $u_2$ coincides with the half-sum of the vector $u_1\bar{u}_2$ and the conjugate vector $u_2\bar{u}_1$:

$$u\bar{u} \equiv \bar{u}u = (u,u)i_0 \quad \Leftrightarrow \quad \frac{u_1\bar{u}_2 + u_2\bar{u}_1}{2} \equiv \frac{\bar{u}_1 u_2 + \bar{u}_2 u_1}{2} = (u_1, u_2)i_0 .$$

From the last relations it is easy to establish the useful identity:

$$u_1\bar{u}u_1 = 2(u_1,u)u_1 - (u_1,u_1)u , \tag{1}$$

where the product $u_1\bar{u}u_1$ is written without brackets, since the resulted vector does not depend on the order of multiplication of hypercomplex factors:

$$u_1\bar{u}u_1 = (u_1\bar{u})u_1 = u_1(\bar{u}u_1).$$

In order to reliably relate to multiplication of quaternions and octonions it is useful to keep in the mind the following elementary formulas.

First formula expresses that conjugation changes the order of factors:

$$\overline{u_1 u_2} = \bar{u}_2\bar{u}_1 .$$

And the following formula for the inner product $(u_1u_2, i_0)$ of vector $u_1u_2$ and vector $i_0$ expresses the transfer rule of the factor from one to another part of the inner product in accompany with simultaneous conjugation:

$$(u_1u_2, i_0) \equiv (u_2u_1, i_0) = (u_1, \bar{u}_2).$$

Finally, the formula:

$$((u_1u)u_2, i_0) = (u_1(uu_2), i_0) \equiv (u_1uu_2, i_0).$$

states that the order of multiplying of three hypercomplex factors $u_1$, $u$ and $u_2$ does not affect the inner product of $(u_1u)u_2$ and $i_0$. So, the product $u_1uu_2$, written without brackets, implies choosing either of the two alternative ways of their arrangement, just as in left part of (1).

To substantiate the claimed additive decomposition in the toolkit of hypercomplex numbers it is useful to preliminary consider a system of three operations of the Hermitian conjugation type.

Let $Au$ be a linear transformation of the vector $u$ of eight-dimensional space of octonions, in particular, the four-dimensional space of quaternions. Let us consider $Au \equiv (u_1\bar{u})u_2$, where $u_1$ and $u_2$ are fixed vectors. Let the operation « $^+$ » of the transformation of the operator $A$ be the operation of Hermitian conjugation:

$$A^+ : \ (Au, v) = (u, A^+v)$$

for any pair of octonions $u$ and $v$.

Let operation «$^*$» be an operation of inversion of multiplicative order, which is carried out by conjugation of $u$ together with both parameters $u_1$ and $u_2$ in combination with common conjugation:

$$A^* : \ A^*u = \overline{(\bar{u}_1\bar{\bar{u}})\bar{u}_2} .$$

Let operation «$^\vee$» be a double conjugation operation, consisting in replacing of the central argument by the conjugate one followed by the common conjugation:

$$A^\vee : \ A^\vee u = \overline{(u_1\bar{\bar{u}})u_2} .$$

The results of the composition of operations are summarized in Tab. 1:

Tab. 1.

An example of three operations forming an Abelian group of mutually inverse elements

| | «» | «$^+$» | «$^*$» | «$^\vee$» |
|---|---|---|---|---|
| «» | $Au \equiv (u_1\bar{u})u_2$ | $A^+u = (u_2\bar{u})u_1$ | $A^*u = u_2(\bar{u}u_1)$ | $A^\vee u = \bar{u}_2(\bar{u}\bar{u}_1)$ |
| «$^+$» | $A^+u = (u_2\bar{u})u_1$ | $A^{++}u = Au$ | $A^{+*}u = u_1(\bar{u}u_2)$ | $A^{+\vee}u = \bar{u}_1(\bar{u}\bar{u}_2)$ |
| «$^*$» | $A^*u = u_2(\bar{u}u_1)$ | $A^{*+}u = u_1(\bar{u}u_2)$ | $A^{**}u = Au$ | $A^{*\vee}u = (\bar{u}_1\bar{u})\bar{u}_2$ |
| «$^\vee$» | $A^\vee u = \bar{u}_2(\bar{u}\bar{u}_1)$ | $A^{\vee+}u = \bar{u}_1(\bar{u}\bar{u}_2)$ | $A^{\vee*}u = (\bar{u}_1\bar{u})\bar{u}_2$ | $A^{\vee\vee}u = Au$ |

Tab. 1 describes the results of application of three operations «$^+$», «$^*$» and «$^\vee$». The table is symmetric and contains the same diagonal elements. This indicates that the operations under consideration form an Abelian group of mutually inverse elements.

## 3. GENERALIZED ADDITIVE DECOMPOSITION OF LINEAR OPERATOR

For several operations of Hermitian conjugation type, which constitute an Abelian group of mutually inverse elements, the linear operator $A$ is trivially decomposed into a sum of mixed *symmetric – skew-symmetric* operators that under the action of each operation either do not change or change the sign.

In a single operation, this decomposition is well known. This is the decomposition of the operator into a Hermitian and a skew-Hermitian parts $(+1)$, $(-1)$:

$$\begin{matrix} A = (+1) + (-1) \\ A^+ = (+1) - (-1) \end{matrix}, \quad (+1), (-1): \quad \begin{matrix} (+1) = \dfrac{A + A^+}{2} \\ (-1) = \dfrac{A - A^+}{2} \end{matrix}$$

For two operations «$^+$» and «$^*$», an analogous decomposition into the sum of *symmetric–skew-symmetric* operators $\begin{pmatrix} \varepsilon_+ = \pm 1 \\ \varepsilon_* = \pm 1 \end{pmatrix}$ is determined by the relations:

$$\begin{pmatrix} +1 \\ +1 \end{pmatrix}, \begin{pmatrix} -1 \\ +1 \end{pmatrix}, \begin{pmatrix} +1 \\ -1 \end{pmatrix}, \begin{pmatrix} -1 \\ -1 \end{pmatrix}: \quad \begin{matrix} A^+\begin{pmatrix} \varepsilon_+ \\ \varepsilon_* \end{pmatrix} = \varepsilon_+\begin{pmatrix} \varepsilon_+ \\ \varepsilon_* \end{pmatrix} \\ A^*\begin{pmatrix} \varepsilon_+ \\ \varepsilon_* \end{pmatrix} = \varepsilon_*\begin{pmatrix} \varepsilon_+ \\ \varepsilon_* \end{pmatrix} \end{matrix}$$

and is written as:

$$\begin{matrix} A \equiv \begin{pmatrix} +1 \\ +1 \end{pmatrix} + \begin{pmatrix} -1 \\ +1 \end{pmatrix} + \begin{pmatrix} +1 \\ -1 \end{pmatrix} + \begin{pmatrix} -1 \\ -1 \end{pmatrix}, & \begin{pmatrix} +1 \\ +1 \end{pmatrix} = \dfrac{A + A^+ + A^* + A^{+*}}{4}, \\ A^+ = \begin{pmatrix} +1 \\ +1 \end{pmatrix} - \begin{pmatrix} -1 \\ +1 \end{pmatrix} + \begin{pmatrix} +1 \\ -1 \end{pmatrix} - \begin{pmatrix} -1 \\ -1 \end{pmatrix}, & \begin{pmatrix} -1 \\ +1 \end{pmatrix} = \dfrac{A - A^+ + A^* - A^{+*}}{4}, \\ A^* = \begin{pmatrix} +1 \\ +1 \end{pmatrix} + \begin{pmatrix} -1 \\ +1 \end{pmatrix} - \begin{pmatrix} +1 \\ -1 \end{pmatrix} - \begin{pmatrix} -1 \\ -1 \end{pmatrix}, & \begin{pmatrix} +1 \\ -1 \end{pmatrix} = \dfrac{A + A^+ - A^* - A^{+*}}{4}, \\ A^{+*} = \begin{pmatrix} +1 \\ +1 \end{pmatrix} - \begin{pmatrix} -1 \\ +1 \end{pmatrix} - \begin{pmatrix} +1 \\ -1 \end{pmatrix} + \begin{pmatrix} -1 \\ -1 \end{pmatrix}, & \begin{pmatrix} -1 \\ -1 \end{pmatrix} = \dfrac{A - A^+ - A^* + A^{+*}}{4}. \end{matrix} \tag{2}$$

The above linear transformation of the operators $A$, $A^+$, $A^*$, $A^\vee$ into symmetric–skew-symmetric operators $\begin{pmatrix} +1 \\ +1 \end{pmatrix}$, $\begin{pmatrix} -1 \\ +1 \end{pmatrix}$, $\begin{pmatrix} +1 \\ -1 \end{pmatrix}$, $\begin{pmatrix} -1 \\ -1 \end{pmatrix}$ and the inverse transformation, are described by the *normalized*[1] symmetric orthogonal Hadamard matrix $A_4 \equiv \begin{bmatrix} +1 & +1 & +1 & +1 \\ +1 & -1 & +1 & -1 \\ +1 & +1 & -1 & -1 \\ +1 & -1 & -1 & +1 \end{bmatrix}$, which coincides, up to a factor, with the inverse matrix: $\frac{A_4^2}{4} = E_4 \equiv \begin{bmatrix} 1 & 0 & 0 & 0 \\ 0 & 1 & 0 & 0 \\ 0 & 0 & 1 & 0 \\ 0 & 0 & 0 & 1 \end{bmatrix}$.

[1] In the *normalized* Hadamard matrix, the first row and the first column consist of only 1.

For a fixed first row in the matrix $A_4$, the permutation of the remaining three rows, in general, violates the symmetry of the matrix, but preserves the set of columns. For example, by interchanging the last two lines, we get a matrix, in which the rows, numbered from top to bottom in order 1, 2, 3, 4, correspond to columns alternating from left to right in the order of 1, 3, 4, 2.

It is characteristic that the quartet of rows (columns) of the matrix $A_4$ is an abelian group with respect to the operation of termwise multiplication. In this case, the rows (columns) 1, $\alpha$, $\beta$, $\alpha\beta$ are alternated in the order of group multiplication which is not violated for any permutation of rows (columns).

For three operations « $^+$ », « $^*$ » and « $^\vee$ », the expansion into a sum of *symmetric – skew-symmetric* operators of the form $\begin{pmatrix} \varepsilon_+ = \pm 1 \\ \varepsilon_* = \pm 1 \\ \varepsilon_\vee = \pm 1 \end{pmatrix}$ whose action is determined by the relations:

$$\begin{pmatrix} +1 \\ +1 \\ +1 \end{pmatrix}, \begin{pmatrix} -1 \\ +1 \\ +1 \end{pmatrix}, \begin{pmatrix} +1 \\ -1 \\ +1 \end{pmatrix}, \begin{pmatrix} -1 \\ -1 \\ +1 \end{pmatrix} \begin{pmatrix} +1 \\ +1 \\ -1 \end{pmatrix}, \begin{pmatrix} -1 \\ +1 \\ -1 \end{pmatrix}, \begin{pmatrix} +1 \\ -1 \\ -1 \end{pmatrix}, \begin{pmatrix} -1 \\ -1 \\ -1 \end{pmatrix} : \quad \begin{matrix} A^+ \begin{pmatrix} \varepsilon_+ \\ \varepsilon_* \\ \varepsilon_\vee \end{pmatrix} = \varepsilon_+ \begin{pmatrix} \varepsilon_+ \\ \varepsilon_* \\ \varepsilon_\vee \end{pmatrix} \\ A^* \begin{pmatrix} \varepsilon_+ \\ \varepsilon_* \\ \varepsilon_\vee \end{pmatrix} = \varepsilon_* \begin{pmatrix} \varepsilon_+ \\ \varepsilon_* \\ \varepsilon_\vee \end{pmatrix} \\ A^\vee \begin{pmatrix} \varepsilon_+ \\ \varepsilon_* \\ \varepsilon_\vee \end{pmatrix} = \varepsilon_\vee \begin{pmatrix} \varepsilon_+ \\ \varepsilon_* \\ \varepsilon_\vee \end{pmatrix} \end{matrix},$$

is presented in the form:

$$A \equiv \begin{pmatrix} +1 \\ +1 \\ +1 \end{pmatrix} + \begin{pmatrix} -1 \\ +1 \\ +1 \end{pmatrix} + \begin{pmatrix} +1 \\ -1 \\ +1 \end{pmatrix} + \begin{pmatrix} -1 \\ -1 \\ +1 \end{pmatrix} + \begin{pmatrix} +1 \\ +1 \\ -1 \end{pmatrix} + \begin{pmatrix} -1 \\ +1 \\ -1 \end{pmatrix} + \begin{pmatrix} +1 \\ -1 \\ -1 \end{pmatrix} + \begin{pmatrix} -1 \\ -1 \\ -1 \end{pmatrix},$$

$$A^+ = \begin{pmatrix} +1 \\ +1 \\ +1 \end{pmatrix} - \begin{pmatrix} -1 \\ +1 \\ +1 \end{pmatrix} + \begin{pmatrix} +1 \\ -1 \\ +1 \end{pmatrix} - \begin{pmatrix} -1 \\ -1 \\ +1 \end{pmatrix} + \begin{pmatrix} +1 \\ +1 \\ -1 \end{pmatrix} - \begin{pmatrix} -1 \\ +1 \\ -1 \end{pmatrix} + \begin{pmatrix} +1 \\ -1 \\ -1 \end{pmatrix} - \begin{pmatrix} -1 \\ -1 \\ -1 \end{pmatrix},$$

$$A^* = \begin{pmatrix} +1 \\ +1 \\ +1 \end{pmatrix} + \begin{pmatrix} -1 \\ +1 \\ +1 \end{pmatrix} - \begin{pmatrix} +1 \\ -1 \\ +1 \end{pmatrix} - \begin{pmatrix} -1 \\ -1 \\ +1 \end{pmatrix} + \begin{pmatrix} +1 \\ +1 \\ -1 \end{pmatrix} + \begin{pmatrix} -1 \\ +1 \\ -1 \end{pmatrix} - \begin{pmatrix} +1 \\ -1 \\ -1 \end{pmatrix} - \begin{pmatrix} -1 \\ -1 \\ -1 \end{pmatrix},$$

$$A^{+*} = \begin{pmatrix} +1 \\ +1 \\ +1 \end{pmatrix} - \begin{pmatrix} -1 \\ +1 \\ +1 \end{pmatrix} - \begin{pmatrix} +1 \\ -1 \\ +1 \end{pmatrix} + \begin{pmatrix} -1 \\ -1 \\ +1 \end{pmatrix} + \begin{pmatrix} +1 \\ +1 \\ -1 \end{pmatrix} - \begin{pmatrix} -1 \\ +1 \\ -1 \end{pmatrix} - \begin{pmatrix} +1 \\ -1 \\ -1 \end{pmatrix} + \begin{pmatrix} -1 \\ -1 \\ -1 \end{pmatrix},$$

$$A^\vee = \begin{pmatrix} +1 \\ +1 \\ +1 \end{pmatrix} + \begin{pmatrix} -1 \\ +1 \\ +1 \end{pmatrix} + \begin{pmatrix} +1 \\ -1 \\ +1 \end{pmatrix} + \begin{pmatrix} -1 \\ -1 \\ +1 \end{pmatrix} - \begin{pmatrix} +1 \\ +1 \\ -1 \end{pmatrix} - \begin{pmatrix} -1 \\ +1 \\ -1 \end{pmatrix} - \begin{pmatrix} +1 \\ -1 \\ -1 \end{pmatrix} - \begin{pmatrix} -1 \\ -1 \\ -1 \end{pmatrix},$$

$$A^{+\vee} = \begin{pmatrix} +1 \\ +1 \\ +1 \end{pmatrix} - \begin{pmatrix} -1 \\ +1 \\ +1 \end{pmatrix} + \begin{pmatrix} +1 \\ -1 \\ +1 \end{pmatrix} - \begin{pmatrix} -1 \\ -1 \\ +1 \end{pmatrix} - \begin{pmatrix} +1 \\ +1 \\ -1 \end{pmatrix} + \begin{pmatrix} -1 \\ +1 \\ -1 \end{pmatrix} - \begin{pmatrix} +1 \\ -1 \\ -1 \end{pmatrix} + \begin{pmatrix} -1 \\ -1 \\ -1 \end{pmatrix},$$

$$A^{*\vee} = \begin{pmatrix}+1\\+1\\+1\end{pmatrix} + \begin{pmatrix}-1\\+1\\+1\end{pmatrix} - \begin{pmatrix}+1\\-1\\+1\end{pmatrix} - \begin{pmatrix}-1\\-1\\+1\end{pmatrix} - \begin{pmatrix}+1\\+1\\-1\end{pmatrix} - \begin{pmatrix}-1\\+1\\-1\end{pmatrix} + \begin{pmatrix}+1\\-1\\-1\end{pmatrix} + \begin{pmatrix}-1\\-1\\-1\end{pmatrix},$$

$$A^{+*\vee} = \begin{pmatrix}+1\\+1\\+1\end{pmatrix} - \begin{pmatrix}-1\\+1\\+1\end{pmatrix} - \begin{pmatrix}+1\\-1\\+1\end{pmatrix} + \begin{pmatrix}-1\\-1\\+1\end{pmatrix} - \begin{pmatrix}+1\\+1\\-1\end{pmatrix} + \begin{pmatrix}-1\\+1\\-1\end{pmatrix} + \begin{pmatrix}+1\\-1\\-1\end{pmatrix} - \begin{pmatrix}-1\\-1\\-1\end{pmatrix}.$$

In this case the matrix $A_8$ of size 8×8, consisting of ± 1, is symmetric and has the form of Fig. 1:

| | $e$ | $\alpha$ | $\beta$ | $\alpha\beta$ | $\gamma$ | $\alpha\gamma$ | $\beta\gamma$ | $\alpha\beta\gamma$ |
|---|---|---|---|---|---|---|---|---|
| $e$ | +1 | +1 | +1 | +1 | +1 | +1 | +1 | +1 |
| $\alpha$ | +1 | -1 | +1 | -1 | +1 | -1 | +1 | -1 |
| $\beta$ | +1 | +1 | -1 | -1 | +1 | +1 | -1 | -1 |
| $\alpha\beta$ | +1 | -1 | -1 | +1 | +1 | -1 | -1 | +1 |
| $\gamma$ | +1 | +1 | +1 | +1 | -1 | -1 | -1 | -1 |
| $\alpha\gamma$ | +1 | -1 | +1 | -1 | -1 | +1 | -1 | +1 |
| $\beta\gamma$ | +1 | +1 | -1 | -1 | -1 | -1 | +1 | +1 |
| $\alpha\beta\gamma$ | +1 | -1 | -1 | +1 | -1 | +1 | +1 | -1 |

Fig. 1. An example of a symmetric matrix with matching sets of rows and columns

Fig. 1 shows a matrix $A_8$ in which the alphabetic designations of the alternation of rows and columns in a fixed order of group multiplication are written before the rows and above the columns. For clarity, the cell fields containing $+1$, are painted black.

A symmetric matrix $A_8$ is a normalized orthogonal Hadamard matrix, and, up to a factor, coincides with its inverse: $A_8^2 / 8 = E_8$, where $E_8$ is the diagonal matrix of 1. The eight rows (columns) in Fig. 1 occupying places $e$, $\alpha$, $\beta$, $\alpha\beta$, $\gamma$, $\alpha\gamma$, $\beta\gamma$, $\alpha\beta\gamma$, constitute an Abelian group with respect to termwise multiplication and are ordered in the matrix accordingly to «doubling algorithm» [2] (Fig. 1).

It is noteworthy that the matrix $A_8$ posesses a *permutation symmetry* [5-9], such that the set of columns does not change under of 168 permutations of rows preserving the order of the alternation of rows in accordance with the doubling algorithm. These 168 permutations constitute a group. Of these, 140 permutations of rows preserve sets of columns but violate the diagonal symmetry of the matrix, as, for example, in the matrix of Fig. 2.

Fig. 2 shows a matrix of identical sets of rows and columns. For clarity, the cell fields of the matrix, containing $+1$, are colored in black. Before the rows and above the columns of the matrix are written the alphabetical designations that indicate the location of rows and columns placed in the order of group term-by-term multiplication. On the left and superiorly the underlined numerical designations of rows and columns are added, which numbered in the order of their alternation in the matrix $A_8$ in Fig. 1. The remaining 28 of the 168 permutations of the rows of the matrix $A_8$ preserve not only the set of columns, but also the diagonal symmetry of the matrix $A_8$ (Fig. 1).

| | | 0 | 3 | 6 | 5 | 7 | 4 | 1 | 2 |
|---|---|---|---|---|---|---|---|---|---|
| | | $e$ | $\alpha$ | $\beta$ | $\alpha\beta$ | $\gamma$ | $\alpha\gamma$ | $\beta\gamma$ | $\alpha\beta\gamma$ |
| 0 | $e$ | +1 | +1 | +1 | +1 | +1 | +1 | +1 | +1 |
| 5 | $\alpha$ | +1 | -1 | +1 | -1 | -1 | +1 | -1 | +1 |
| 7 | $\beta$ | +1 | -1 | -1 | +1 | -1 | +1 | +1 | -1 |
| 2 | $\alpha\beta$ | +1 | +1 | -1 | -1 | +1 | +1 | -1 | -1 |
| 6 | $\gamma$ | +1 | +1 | -1 | -1 | -1 | -1 | +1 | +1 |
| 3 | $\alpha\gamma$ | +1 | -1 | -1 | +1 | +1 | -1 | -1 | +1 |
| 1 | $\beta\gamma$ | +1 | -1 | +1 | -1 | +1 | -1 | +1 | -1 |
| 4 | $\alpha\beta\gamma$ | +1 | +1 | +1 | +1 | -1 | -1 | -1 | -1 |

Fig. 2. Example of an asymmetric matrix with matching rows and columns

The alternation of rows 0,1, 2, 3, 4, 5, 6, 7 (Fig. 1, 2) in symmetric matrices is given in Tab. 2.

Tab. 2.

Alternating of rows in 28 symmetric matrices

| | | | | | | | | | | | | | | | | | | | | | | | | | | | | |
|---|---|---|---|---|---|---|---|---|---|---|---|---|---|---|---|---|---|---|---|---|---|---|---|---|---|---|---|---|
| 0 | 0 | 0 | 0 | 0 | 0 | 0 | 0 | 0 | 0 | 0 | 0 | 0 | 0 | 0 | 0 | 0 | 0 | 0 | 0 | 0 | 0 | 0 | 0 | 0 | 0 | 0 | 0 | 0 |
| 1 | 1 | 1 | 2 | 3 | 4 | 5 | 3 | 4 | 6 | 5 | 2 | 6 | 7 | 1 | 7 | 2 | 7 | 3 | 5 | 4 | 6 | 2 | 5 | 4 | 6 | 3 | 7 |
| 2 | 4 | 6 | 3 | 1 | 2 | 2 | 7 | 2 | 7 | 6 | 1 | 5 | 3 | 4 | 5 | 7 | 3 | 5 | 4 | 6 | 1 | 5 | 4 | 6 | 3 | 7 | 1 |
| 3 | 5 | 7 | 1 | 2 | 6 | 7 | 4 | 6 | 1 | 3 | 3 | 3 | 4 | 5 | 2 | 5 | 4 | 6 | 1 | 2 | 7 | 7 | 1 | 2 | 5 | 4 | 6 |
| 4 | 6 | 2 | 4 | 4 | 5 | 1 | 6 | 1 | 3 | 7 | 4 | 7 | 5 | 2 | 3 | 6 | 1 | 2 | 7 | 3 | 5 | 6 | 3 | 7 | 1 | 2 | 5 |
| 5 | 7 | 3 | 6 | 7 | 1 | 4 | 5 | 5 | 5 | 2 | 6 | 1 | 2 | 3 | 4 | 4 | 6 | 1 | 2 | 7 | 3 | 4 | 6 | 3 | 7 | 1 | 2 |
| 6 | 2 | 4 | 7 | 5 | 7 | 3 | 1 | 3 | 4 | 1 | 5 | 2 | 6 | 6 | 6 | 1 | 2 | 7 | 3 | 5 | 4 | 3 | 7 | 1 | 2 | 5 | 4 |
| 7 | 3 | 5 | 5 | 6 | 3 | 6 | 2 | 7 | 2 | 4 | 7 | 4 | 1 | 7 | 1 | 3 | 5 | 4 | 6 | 1 | 2 | 1 | 2 | 5 | 4 | 6 | 3 |

In the above table the leftmost column of numbers, alternating from top to bottom in a natural order, describes the matrix of Fig. 1. Pairs of columns from Tab. 2 define the permutations of the rows preserving the symmetry of matrix Fig.1. For a fixed first column that coincides with a column of alternating numbers 0,1, 2, 3, 4, 5, 6, 7, the columns of Tab. 2 describe 28 permutations of rows preserving the diagonal symmetry of the matrix in Fig. 1. These 28 permutations are not a group, but represent the union of an octet of cyclic subgroups of 3, 4 and 7 order, including two subgroups of the seventh order, three cyclic subgroups of the third order and three cyclic subgroups of the fourth order. Vertical lines in Tab. 2 separate columns that describe permutations from one cyclic subgroup.

It can be shown that the discussing normalized Hadamard matrices of size 8×8:

- are determined by their *permutability* so that the matrix itself is uniquely reconstructed from permutations that do not violate its diagonal symmetry, and to recover the matrix, it is sufficient to select in prescribed algorithm from twenty-eight columns of Tab. 2 only three columns that specify a pair of permutations [6-8];
- are characterized by the maximum permutation symmetry, which is preserved with the maximum number of permutations of the rows different from the first row containing only +1 [7].

Owing to listed properties, the discussed Hadamard matrices have found practical application in a series of games with permutative symmetry [6-9].

## 4. THE DECOMPOSITION OF THE TRIPLE OCTONIONIC PRODUCT INTO THE SUM OF THE ANTICOMMUTATOR, COMMUTATOR, AND ASSOCIATOR

Let's obtain an additive decomposition of the product $(u_1\bar{u})u_2$ of three octonions $u_1$, $u$ and $u_2$ for two commuting operations, namely, for the operation of Hermitian conjugation « $^+$ » and the operation of inversion of the multiplicative order « $^*$ ».

Using Tab.1, we rewrite the trivial relations from the right-hand side of (2) as:

$$\begin{pmatrix}+1\\+1\end{pmatrix}u \equiv \{u_1,u,u_2\} = \frac{(u_1\bar{u})u_2 + (u_2\bar{u})u_1}{2} = \frac{u_1(\bar{u}u_2) + u_2(\bar{u}u_1)}{2}, \tag{3}$$

$$\begin{pmatrix}-1\\+1\end{pmatrix}u \equiv \langle u_1,u,u_2\rangle = \frac{(u_1\bar{u})u_2 - u_1(\bar{u}u_2)}{2} = \frac{-(u_2\bar{u})u_1 + u_2(\bar{u}u_1)}{2}, \tag{4}$$

$$\begin{pmatrix}+1\\-1\end{pmatrix}u = 0 = \frac{(u_1\bar{u})u_2 - u_1(\bar{u}u_2)}{4} + \frac{(u_2\bar{u})u_1 - u_2(\bar{u}u_1)}{4} \equiv \langle u_1,u,u_2\rangle + \langle u_2,u,u_1\rangle, \tag{5}$$

$$\begin{pmatrix}-1\\-1\end{pmatrix}u \equiv [u_1,u,u_2] = \frac{(u_1\bar{u})u_2 - u_2(\bar{u}u_1)}{2} = \frac{u_1(\bar{u}u_2) - (u_2\bar{u})u_1}{2}, \tag{6}$$

where the term $\begin{pmatrix}+1\\-1\end{pmatrix}u$ vanishes, and the remaining terms of the additive decomposition receive the names that reflect their expression by commutating the arguments $u_1$, $u_2$, and changing by the opposite the multiplicative order of $u_1$, $u$ and $u_2$. So, matching to the traditional terminology, $\{u_1,u,u_2\}$ is called a commutator, $\langle u_1,u,u_2\rangle$ is called an associator, and $[u_1,u,u_2]$ is called an anticommutator. According to (3), (4), (6), they coincide with the symmetric–skew-symmetric (antisymmetric) components $\begin{pmatrix}+1\\+1\end{pmatrix}u$, $\begin{pmatrix}-1\\+1\end{pmatrix}u$, $\begin{pmatrix}-1\\-1\end{pmatrix}u$ of the additive decomposition of the triple product $(u_1\bar{u})u_2$ of octonions $u_1$, $u$ and $u_2$ with the conjugate central argument $u$.

NOTE: The validity of (5) and the equalities of the half-sums in (3), (4), (6) follows from the anticommutativity of the associator, which in turn trivially follows from associator zeroing if any two of the three arguments coincide [10-12]. Concerning formulas (6), it should be noted that they can be considered as a definition of a triple commutator. It is noteworthy that similar definitions are found in the literature, for example, in [13] on the page 224.

By means of the inverse expressions from the left-hand side of (2) the terms $Au = (u_1\bar{u})u_2$, $A^+u = (u_2\bar{u})u_1$, $A^*u = u_2(\bar{u}u_1)$ and $A^{+^*}u = u_1(\bar{u}u_2)$ (see Tab. 1) will be written via $\begin{pmatrix}+1\\+1\end{pmatrix}u \equiv \{u_1,u,u_2\}$, $\begin{pmatrix}-1\\+1\end{pmatrix}u \equiv \langle u_1,u,u_2\rangle$, $\begin{pmatrix}+1\\-1\end{pmatrix}u = 0$, $\begin{pmatrix}-1\\-1\end{pmatrix}u \equiv [u_1,u,u_2]$ as:

$$(u_1\bar{u})u_2 \equiv \{u_1,u,u_2\} + [u_1,u,u_2] + \langle u_1,u,u_2\rangle, \tag{7}$$

$$(u_2\bar{u})u_1 = \{u_1,u,u_2\} - [u_1,u,u_2] - \langle u_1,u,u_2\rangle, \tag{8}$$

$$u_2(\bar{u}u_1) = \{u_1,u,u_2\} - [u_1,u,u_2] + \langle u_1,u,u_2\rangle, \tag{9}$$

$$u_1(\bar{u}u_2) = \{u_1,u,u_2\} + [u_1,u,u_2] - \langle u_1,u,u_2\rangle. \tag{10}$$

Formula (7) gives the decomposition of the triple octonionic product $(u_1\bar{u})u_2$ into the sum of three mutually orthogonal terms, namely, anticommutator $\{u_1,u,u_2\}$, commutator $[u_1,u,u_2]$ and associator $\langle u_1,u,u_2\rangle$:

$$(\{u_1,u,u_2\},[u_1,u,u_2]) = (\{u_1,u,u_2\},\langle u_1,u,u_2\rangle) = ([u_1,u,u_2],\langle u_1,u,u_2\rangle) = 0,$$

and (8) – (10) give the expressions $(u_2\bar{u})u_1$, $u_2(\bar{u}u_1)$ и $u_1(\bar{u}u_2)$ as a linear combinations of the listed mutually orthogonal terms.

It is easy to verify the mutual orthogonality of $\{u_1,u,u_2\}$, $[u_1,u,u_2]$ и $\langle u_1,u,u_2\rangle$, using the explicit expression of $\{u_1,u,u_3\}$ via $u_1$, $u$ and $u_2$:

$$\{u_1,u,u_2\}=(u_1,u)u_2-(u_1,u_2)u+(u,u_2)u_1, \tag{11}$$

and also by the expression of $[u_1,u,u_2]$ via of the multiplicative unit $i_0$ and the cross products $[u,u_2]$, $[u_1,u_2]$ and $[u_1,u]$:

$$[u_1,u,u_2]=([u_1,u],u_2)i_0-(u_1,i_0)[u,u_2]+(u,i_0)[u_1,u_2]-(u_2,i_0)[u_1,u]. \tag{12}$$

The formula (11) for $\{u_1,u,u_2\}$ is derived from identity (1) by replacing of $u_1$ with $u_1+u_2$. The formula (12) for $[u_1,u,u_2]$ is derived from its expression (6) as one or another of the half-sum and the following expression of the product of two hypercomplex numbers via $[u_1,u_2]$ [14]:

$$u_1u_2=(u_1,i_0)u_2+(u_2,i_0)u_1-(u_1,u_2)i_0+[u_1,u_2]. \tag{14}$$

The cross product $[u_1,u_2]$ of two arguments is obtained by replacing in $[u_1,u,u_2]$ the central argument $u$ by the multiplicative unit $i_0$, that according to (6) coincides with the commutator of two vectors $u_1$, и $u_2$:

$$[u_1,u_2]\equiv[u_1,i_0,u_2]=\frac{u_1u_2-u_2u_1}{2}, \tag{13}$$

where the vector $[u_1,u_2]$ is the cross product of a pair of hypercomplex vectors $u_1$ and $u_2$ [2, 4].

The cross product $[u_1,u_2]$ of the hypercomplex numbers $u_1$ and $u_2$ is orthogonal to the multiplicative unity $i_0$ and vanishes when $i_0$ is taken as one or another argument:

$$([u_1,u_2],i_0)=0\ ,$$
$$[i_0,u_2]=[u_1,i_0]=0\ .$$

In the rest, the cross product $[u_1,u_2]$ preserves the properties of the conventional cross product, which is introduced in three-dimensional space using the intuitively perceived "right-hand rule". Under the preservation of properties it is meant that if we introduce the operation of annulling the real component of a hypercomplex number $u$ denoting by a prime:

$$u'=u-(u,i_0)i_0,$$

then the cross product $[u_1,u_2]$ of the vectors $u_1$ and $u_2$ coincides with the conventional cross product $[u_1',u_2']$ of vectors $u_1'$ and $u_2'$ of the three-dimensional subspace that is orthogonal to a multiplicative unit $i_0$.

Using (12), it is easy to establish the orthogonality of the $[u_1,u,u_2]$ to each of three arguments $u_1$, $u$ and $u_2$:

$$([u_1,u,u_2],u_1)=([u_1,u,u_2],u)=([u_1,u,u_2],u_2)=0$$

and to verify that $[u_1,u,u_2]$ and $\{u_1,u,u_2\}$ are mutually orthogonal using (11). The orthogonality of the associator $\langle u_1,u,u_2\rangle$ to each of the three arguments $u_1$, $u$ and $u_2$ is derived from its expression by means of half-sums (4). The same property of the orthogonality of the associator $\langle u_1,u,u_2\rangle$ to each of its arguments $u_1$, $u$ and $u_2$:

$$(\langle u_1,u,u_2\rangle,u_1)=(\langle u_1,u,u_2\rangle,u)=(\langle u_1,u,u_2\rangle,u_2)=0$$

is deduced from its expression using the half-sum (4).

The orthogonality of the triple commutator $[u_1,u,u_2]$ and the associator $\langle u_1,u,u_2\rangle$ to its arguments implies the anticommutativity of mixed products, for example:

$$([u_1,u,u_2],u_3)=-([u_3,u,u_2],u_1),$$
$$(\langle u_1,u,u_2\rangle,u_3)=-(\langle u_3,u,u_2\rangle,u_1).$$

From (12) it is easy to establish that the triple vector product $[u_1,u,u_2]$ vanishes when a pair of its arguments coincide with one another, which implies its anticommutativity, for example:

$$[u_1,u,u_1]=0 \quad \Leftrightarrow \quad [u_1,u,u_2]=-[u_2,u,u_1].$$

The associator $\langle u_1,u,u_2\rangle$ has the same properties of vanishing to zero and anticommutativity, for example:

$$\langle u_1,u,u_1\rangle=0 \quad \Leftrightarrow \quad \langle u_1,u,u_2\rangle=-\langle u_2,u,u_1\rangle.$$

Moreover, it is known [10-12] that the associator vanishes in the case of three arguments from one quaternionic subalgebra. Then it follows from the anticommutativity of the mixed product that the associator $\langle u_1,u,u_2\rangle$ is orthogonal to the vectors of quaternionic subalgebras generated by any pair of three vectors $u_1$, $u$ and $u_2$. Therefore, in addition to the vectors $u_1$, $u$ and $u_2$ the associator $\langle u_1,u,u_2\rangle$ is also orthogonal to the four following vectors:

$$(\langle u_1,u,u_2\rangle,i_0)=(\langle u_1,u,u_2\rangle,[u_1,u])=(\langle u_1,u,u_2\rangle,[u_1,u_2])=(\langle u_1,u,u_2\rangle,[u,u_2])=0.$$

In particular, the associator $\langle u_1,u,u_2\rangle$ is orthogonal to the triple anticommutator $\{u_1,u,u_2\}$ and commutator $[u_1,u,u_2]$, which are linear combinations (11), (12) of the specified seven vectors.

For brevity, further formulas are given without proof, since they are obtained by applying the previously written ones.

NOTE: In calculations, it is recommended to avoid the rut utilizing of the orts $i$, $j$, $k$ etc. [1, 12] without a special need for a spatial basis, since this leads to cumbersome expressions, especially when working with octonions.

## 5. INTERPRETATION

Interpretation of terms $\{u_1,u,u_2\}$, $[u_1,u,u_2]$ и $\langle u_1,u,u_2\rangle$ in the additive decomposition (7) - (10) of products of three octonions $u_1$, $u$ and $u_2$ with conjugate central argument $u$ is as follows.

According to the defining relations (3), (11) $\{u_1,u,u_2\}$ is a *triple anticommutator*, which is expressed as a linear combination of the arguments $u_1$, $u$, $u_2$ and does not change when the first and last arguments are interchanged. When the central argument is replaced by the multiplicative unit $i_0$, the anticommutator $\{u_1,u,u_2\}$ is converted into an anticommutator $\{u_1,u_2\}$ for the product of two hypercomplex numbers [2,14].

The *triple commutator* $[u_1,u,u_2]$ introduced in (5) and represented in (12) as a linear combination of the multiplicative unit $i_0$ and pairwise cross products of the arguments is treated as a generalized cross product of three hypercomplex numbers (quaternions or octonions). When one of the arguments is replaced by the multiplicative unit $i_0$, the triple cross product up to a sign coincides with the conventional cross product of two other arguments. When the central argument $u$ is substituted by the multiplicative unit $i_0$, the triple commutator $[u_1,u,u_2]$ is converted into an ordinary cross product $[u_1,u_2]$ of two arguments $u_1$ and $u_2$ [2].

The associator $\langle u_1,u,u_2\rangle$ is well known in the scientific papers [10-12]. Here it is introduced by formulas (4), where, unlike the conventional definitions, conjugation of the central argument is used, and the associator $\langle u_1,u,u_2\rangle$ is calculated with the coefficient $1/2$. These differences do not affect the basic properties of the associator.

## 6. LENGTHS OF ADDITIVE DECOMPOSITION TERMS

The length squares of the triple anticommutator $\{u_1,u,u_2\}$, the commutator (triple cross product) $[u_1,u,u_2]$, and the associator $\langle u_1,u,u_2\rangle$ are expressed as follows:

$$\|\{u_1,u,u_2\}\|^2 = (u_1,u_1)(u,u)(u_2,u_2) - \det\begin{bmatrix}(u_1,u_1) & (u_1,u) & (u_1,u_2)\\ (u_1,u) & (u,u) & (u,u_2)\\ (u_1,u_2) & (u,u_2) & (u_2,u_2)\end{bmatrix}, \tag{15}$$

$$\|[u_1,u,u_2]\|^2 = ([u_1,u],u_2)^2 + \det\begin{bmatrix}(u_1,u_1) & (u_1,u) & (u_1,u_2)\\ (u_1,u) & (u,u) & (u,u_2)\\ (u_1,u_2) & (u,u_3) & (u_2,u_2)\end{bmatrix} - \det\begin{bmatrix}(u_1',u_1') & (u_1',u') & (u_1',u_2')\\ (u_1',u') & (u',u') & (u',u_2')\\ (u_1',u_2') & (u',u_2') & (u_2',u_2')\end{bmatrix}, \tag{16}$$

$$\|\langle u_1,u,u_2\rangle\|^2 = \det\begin{bmatrix}(u_1',u_1') & (u_1',u') & (u_1',u_2')\\ (u_1',u') & (u',u') & (u',u_2')\\ (u_1',u_2') & (u',u_2') & (u_2',u_2')\end{bmatrix} - ([u_1,u],u_2)^2. \tag{17}$$

where for any hypercomplex vector $u$ the designation $\|u\|$ denotes its length, and $\|u\|^2$ is used as an equivalent of $(u,u)$ to shorten the record. So in (15)-(17):

$$\|\{u_1,u,u_2\}\|^2 \equiv (\{u_1,u,u_2\},\{u_1,u,u_2\}),$$

$$\|[u_1,u,u_2]\|^2 \equiv ([u_1,u,u_2],[u_1,u,u_2]),$$

$$\|\langle u_1,u,u_2\rangle\|^2 \equiv (\langle u_1,u,u_2\rangle,\langle u_1,u,u_2\rangle).$$

It is useful to note that the length square of the anticommutative component $[u_1,u,u_2]+\langle u_1,u,u_2\rangle$ of the triple octonionic product $(u_1\bar{u})u_2$ is expressed by the determinant of the symmetric positive definite matrix of scalar products of the arguments $u_1$, $u$ and $u_2$:

$$\|[u_1,u,u_2]+\langle u_1,u,u_2\rangle\|^2 = \det\begin{bmatrix}(u_1,u_1) & (u_1,u) & (u_1,u_2)\\ (u_1,u) & (u,u) & (u,u_2)\\ (u_1,u_2) & (u,u_2) & (u_3,u_2)\end{bmatrix}.$$

NOTE: In applications, it may be useful to represent (16), (17) in another form, using an expression of $\det\begin{bmatrix}(u_1',u_1') & (u_1',u_2') & (u_1',u_3')\\ (u_1',u_2') & (u_2',u_2') & (u_2',u_3')\\ (u_1',u_3') & (u_2',u_3') & (u_3',u_3')\end{bmatrix}$ in terms of pairwise scalar products and scalar products of one argument to the conjugate second argument:

$$\det\begin{bmatrix}(u_1',u_1') & (u_1',u') & (u_1',u_2')\\ (u_1',u') & (u',u') & (u',u_2')\\ (u_1',u_2') & (u',u_2') & (u_2',u_2')\end{bmatrix} = \frac{\det\begin{bmatrix}(u_1,u_1) & (u_1,u) & (u_1,u_2)\\ (u_1,u) & (u,u) & (u,u_2)\\ (u_1,u_2) & (u,u_2) & (u_2,u_2)\end{bmatrix} - \det\begin{bmatrix}(u_1,\bar{u}_1) & (u_1,\bar{u}) & (u_1,\bar{u}_2)\\ (u_1,\bar{u}) & (u,\bar{u}) & (u,\bar{u}_2)\\ (u_1,\bar{u}_2) & (u,\bar{u}_2) & (u_2,\bar{u}_2)\end{bmatrix}}{2}.$$

It is pertinent to note that in the special theory of relativity of A. Einstein, the quantity $(u,\bar{u})$ is called «the space-time interval» and is treated as an analog of the spatial distance, where $u$ is the four-dimensional space-time vector identified with the hypercomplex number by W.R. Hamilton himself and by his contemporary followers [1,15].

From the decomposition (7) of the product $(u_1\bar{u})u_3$ into orthogonal terms by means of the summation (15)-(17), we obtain the equality:

$$\|(u_1\bar{u})u_2\|^2 \equiv ((u_1\bar{u})u_2,(u_1\bar{u})u_2) = (u_1,u_1)(u,u)(u_2,u_2),$$

expressing the axiomatic property of *normed* algebras, which consists in the fact that the square of the product length of two hypercomplex vectors is equal to the product of the squares of the lengths of these vectors: $(u_1u_2,u_1u_2) = (u_1,u_1)(u_2,u_2) \Leftrightarrow ((u_1\bar{u})u_2,(u_1\bar{u})u_2) = (u_1,u_1)(u,u)(u_2,u_2)$.

## 7. SOURCE SOLUTIONS

The generalizing of cross product of two vectors to the case of more than three dimensions was investigated by Zurab K. Silagadze in [16]. The solution is constructed by postulating axioms, which have a physical meaning and generalize the characteristic properties of the traditional three-dimensional cross product. It is shown that the only possible dimensionality for a justified generalization of a vector product is dimension 7. The generalization of a cross product as the octonionic commutator of two arguments is given as an example. It is noteworthy that in addition to the product of three vectors in [16], the product of three arguments is also considered. In this connection, a number of interesting relations are derived. Namely the known formula "'BAC' minus 'CAB'" (russian) is treated:

$$[A,[B,C]]-B(A,C)+C(A,B)\equiv\langle A,B,C\rangle,$$

where, in the notation used here $A\equiv u_1'$, $B\equiv u'$ and $C\equiv u_2'$

For comparison of formulas, the difference in notations is presented in Tab. 3.

Tab. 3.

Matching of notations [16]

| Notion | Notation | |
|---|---|---|
| | In [16] | In present paper |
| Multiplicative unit | $e$ | $i_0$ |
| Inner product ("Euclidean scalar product") | $\vec{u}_1\cdot\vec{u}_2$ | $(u_1',u_2')$ |
| Cross product of two arguments ("Vector product") | $\vec{u}_1\times\vec{u}_2$ | $[u_1',u_2']$ |
| Associator ("Ternary product") | $-\{\vec{u}_1,\vec{u},\vec{u}_2\}$ | $\langle u_1',u',u_2'\rangle$ |

In [4] Susumu Okubo considered the product $(u_1u)u_2$ of three octonions $u_1$, $u$, $u_2$, extracting out the anticommutative part $[u_1,u,u_2]_{Ocubo}$ in the additive decomposition. The solution obtained by a fairly long chain of calculations is given on the page 22 in the form:

$$(u_1u)u_2=[u_1,u,u_2]_{Ocubo}+2(u,i_0)u_1u_2-(u,u_2)u_1-(u_1,u)u_2+(u_1,u_2)u,$$

wherein $[u_1,u,u_2]_{Okubo}$ is disclosed in our notations by the following expression:

$$[u_1,u,u_2]_{Okubo}=-\langle u_1,u,u_2\rangle+(u_1,i_0)[u,u_2]+(u,i_0)[u_2,u_1]+(u_2,i_0)[u_1,u]-(u_2,[u_1,u])i_0.$$

To simplify the comparison of the decomposition of the triple octonionic product $(u_1\bar{u})u_2$ in the present paper with the decomposition of the triple octonionic product $(u_1u)u_2$ in [4] we give the Tab. 4 of co-assignments of notations.

Tab. 4.

Matching of notations [4]

| Notion | Notation | |
|---|---|---|
| | In [4] | In present paper |
| Multiplicative unit | $e$ | $i_0$ |
| Inner product | $\langle u_1\vert u_2\rangle$ | $(u_1,u_2)$ |
| Cross product of two arguments | $\frac{[u_1,u_2]}{2}$ | $[u_1,u_2]$ |
| Associator | $-\frac{(u_1,u,u_2)}{2}$ | $\langle u_1,u,u_2\rangle$ |

Taking into account the expressions (11) for $\{u_1,u,u_2\}$ and (12) for $[u_1,u,u_2]$, the solution [4] is expressed by the formula:

$$(u_1u)u_2=2(u,i_0)u_1u_2-\{u_1,u,u_2\}-[u_1,u,u_2]-\langle u_1,u,u_2\rangle,$$

which is equivalent to the formula (7) for the decomposition of the triple product $(u_1\bar{u})u_2$ into the sum of triple anticommutator $\{u_1,u,u_2\}$, triple commutator $[u_1,u,u_2]$, and associator $\langle u_1,u,u_2\rangle$.

Thus, the formulas (7)-(10) improve the results, originally obtained in [4] for $(u_1u)u_2$ and allow to trivially express the idea of decomposition of $(u_1\bar{u})u_2$ into the sum of mutually orthogonal $\{u_1,u,u_2\}$, $[u_1,u,u_2]$ and $\langle u_1,u,u_2\rangle$ by commutating of factors and inversing of the multiplicative order of the arguments $u_1$, $u$ and $u_2$. At the same time, the identities (3)-(6) develop the interpretation of $\{u_1,u,u_2\}$, $[u_1,u,u_2]$ and $\langle u_1,u,u_2\rangle$ as symmetric–antisymmetric or, avoiding ambiguous terminology, *symmetric–skew-symmetric* parts of the product $(u_1\bar{u})u_2$.

In contrast to [4], in [17] the authors Tevian Dray and Corinne A. Manogue consider the product $u_1(\bar{u}u_2)$ of octonions $u_1$, $u$ and $u_2$ with the conjugate central element $\bar{u}$, but the triple octonion product denoted $u_1\times u\times u_2$ they defined as:

$$u_1\times u\times u_2\equiv\frac{u_1(\bar{u}u_2)-u_2(\bar{u}u_1)}{2}=[u_1,u,u_2]-\langle u_1,u,u_2\rangle,$$

where the generalized cross product $[u_1,u,u_2]$ and the associator $\langle u_1,u,u_2\rangle$ are combined in one difference expression that obtained from (9), (10). Such the product $u_1\times u\times u_2$ is anticommutative, but when the product $u_1(\bar{u}u_2)$ is decomposed into a sum of three orthogonal terms, it does not lead to a transparent expression (7), justified by trivial commutation rules and rules for inversing of the multiplicative order of the arguments (3)-(6).

## 8. SYMMETRY CONSIDERATIONS

The paper [16] is well known in abstract algebra (section RA – Rings and Algebras). In this paper, we draw the attention of the reader to the remarkable results in this field of the well-known physics Susum Okubo, who in 1993, in a triple octonionic product has extracted out two additive anticommutative parts guiding by symmetry considerations. In the continuation of S. Okubo's research, one of the parts is proposed to be interpreted as a generalization of the cross product for four-dimensional and eight-dimensional Euclidean space and the case of three arguments basing on symmetry considerations. It is important that the notion of symmetry is associated with the alternative notion of antisymmetry and is interpreted in the following senses:

- as commutativity and anticommutativity of the ternary operation with vector arguments, in which the permutation of the arguments either does not affect the result of the operation, or changes the sign of the resultant vector to the opposite;
- as the properties of several commuting operations on linear operators in a vector space, which are introduced as a generalization of the Hermitian decomposition of the operator into a symmetric and antisymmetric (skew-symmetric) parts.

So, in this paper an additive decomposition of the product of three octonions with conjugate central argument is obtained. The main results are as follows:

- The additive decomposition consists of three parts, namely, the triple anticommutator $\{u_1,u,u_2\}$, the commutator (triple vector product) $[u_1,u,u_2]$, and the associator $\langle u_1,u,u_2\rangle$, determined by changing the sequence order and the order of multiplication of the three arguments $u_1, u$, $u_2$;
- all three terms of the additive decomposition are mutually orthogonal;
- two terms: $[u_1,u,u_2]$ and $\langle u_1,u,u_2\rangle$ are anticommutative;
- all three terms (with zero fourth term $\begin{pmatrix}+1\\-1\end{pmatrix}u=0$) are symmetric–skew-symmetric components $\{u_1,u,u_2\}=\begin{pmatrix}+1\\+1\end{pmatrix}u$, $[u_1,u,u_2]=\begin{pmatrix}-1\\-1\end{pmatrix}u$, $\langle u_1,u,u_2\rangle=\begin{pmatrix}-1\\+1\end{pmatrix}u$ of the representations of the operator of multiplication of two constant octonions $u_1$, $u_2$ by a central conjugate octonion $\bar{u}$, considered as a variable argument, where the upper and lower signs inside the brackets denote the property of mixed symmetry either to preserve, or change the sign under the operation of Hermitian conjugation and either preserve or change the sign under the operation of inversion of the multiplicative order, respectively.

The additional symmetry consideration, obtained and patented as a by-product [7, 8], is that the symmetric normalized Hadamard matrices, arising in the generalization of the additive Hermitian decomposition of an operator into symmetric–skew-symmetric parts, have the property of preserving the mirror diagonal symmetry under the maximum number of row or column permutations.

Symmetry considerations provide a proper choice of the definition of a generalized triple vector product from the number of possible generalizations. Thus, the term $[u_1, u, u_2]$, most likely, is the appropriate variant of the triple cross product, generalized to the case of four-dimensional and eight-dimensional Euclidean space.

## 9. CONCLUSION

In addition to the cross product definition, quaternions and octonions are famous for the fact that they gracefully describe the rotation of space orthogonal to the multiplicative unit $i_0$. In a later paper we are going to show that this also applies to the Lorentz transformations, which in this case are considered in Euclidean space [18-20]. Apparently, the refinement of the representation of Lorentz transformations in terms of quaternions and octonions is promising for application in Einstein's theory of relativity [15], Dirac quantum mechanics [15,21] and other fields of physics [4]. However, basing on the common experience of our preceding work, the development of applications in the field of digital image processing is closer to us. Therefore, it is planned, first of all, to continue research in this context, relying on our own experience in image processing, as well as on long-term studies of the Yoshkar-Ola's school of scientists who successfully apply quaternions for image processing [22-24].

As is known in image processing, the pixel of the color image is mapped to the point of the three-dimensional Euclidean color space presented in one or another coordinate system. To describe the color transformation of an image, the rotation of the color space is usually used. Using the Lorentz transformation instead of simple rotation, we will construct, programmatically implement and explore the generalized color transformation of the image to advance in the practical solution of one of the many problems of modern image processing.